\documentclass[12pt]{amsart}
\usepackage{latexsym,amscd,amssymb,epsfig,enumerate} 
\usepackage[dvips]{color}

\theoremstyle{plain}
  \newtheorem{theorem}{Theorem}[section]
  \newtheorem{proposition}[theorem]{Proposition}
  \newtheorem{lemma}[theorem]{Lemma}
  \newtheorem{corollary}[theorem]{Corollary}
  \newtheorem{conjecture}[theorem]{Conjecture}
\theoremstyle{definition}
  \newtheorem{definition}[theorem]{Definition}
  \newtheorem{example}[theorem]{Example}
  
  \newtheorem{question}[theorem]{Question}
 \theoremstyle{remark}
  \newtheorem{remark}[theorem]{Remark}

\numberwithin{equation}{section}

\def\EOP{\hfill$\Box$}

\def\integers{{\mathbb Z}}

\def\mm{{\mathbf m}}

\def\capac{cap}

\begin{document}

\title{Shelling Coxeter-like complexes and sorting on trees}

\author{Patricia Hersh}
\address{Department of Mathematics\\
     Indiana University\\
     Rawles Hall\\
     Bloomington, IN 47405}
\email{phersh@indiana.edu}
\address{Box 8205\\
Department of Mathematics\\
North Carolina State University\\
Raleigh, NC 27695-8205}

\thanks{The author was supported by NSF grants  DMS-0500638 and DMS-0757935.}

\begin{abstract}
In their work on `Coxeter-like complexes', Babson and Reiner
introduced a simplicial complex $\Delta_T$
associated to each tree $T$ on $n$ nodes, 
generalizing chessboard complexes and type A Coxeter complexes.
They conjectured that $\Delta_T$ is $(n-b-1)$-connected when the 
tree has $b$ leaves.  We provide a shelling for 
the $(n-b)$-skeleton of $\Delta_T$, thereby proving this 
conjecture.

In the process, we introduce notions of weak order and inversion functions 
on the labellings of a tree $T$ which imply shellability of $\Delta_T$, and  we
construct such inversion functions for a large enough class of trees to deduce
the aforementioned conjecture and also recover the shellability of 
chessboard complexes $M_{m,n}$ with $n \ge 2m-1$.  We also prove that the existence or
nonexistence of an inversion function for a fixed tree governs which networks with a tree
structure admit greedy sorting algorithms by inversion elimination and
provide an inversion function for trees where 
each vertex has capacity at least its degree minus one.
\end{abstract}

\maketitle

\section{Introduction.}

Coxeter-like complexes were introduced in [BR]
as a common generalization of Coxeter complexes and chessboard complexes.
The point was to associate a cell complex analogous to a Coxeter 
complex to any minimal generating set for any finite group 
-- for instance,
to any set of $n-1$ transpositions in the symmetric group $S_n$ which
generate $S_n$.   While traditional Coxeter complexes are endowed with a wealth
of beautiful and remarkable properties (see e.g. [Hu],  [Br]), including 
shellability (as proven 
in [Bj2]), the topological structure of more general 
Coxeter-like complexes is often much more subtle.
The chessboard complexes already demonstrate this.

Chessboard complexes have been studied extensively, motivated both by
applications to computational geometry (see e.g. [ZV]) and also 
because of a direct relation established in [RR] between their 
homology groups and the
Tor groups of  Segre modules.
The homology groups for chessboard complexes are not easy to 
determine, for example involving 3-torsion; results in [BLVZ] and [SW] together show
that they are exactly $(\nu_{m,n}-1)$-connected, for 
$\nu_{m,n} = \min\{m,n,\lfloor \frac{m+n+1}{3}\rfloor \} -1$.   Our  main focus will be on 
topological properties of more general Coxeter-like complexes.
Let us now briefly recall some 
terminology and establish some notation, before describing our  main results.

Babson and Reiner associate to any finite group $G$ with chosen
minimal generating set
$S$ a cell complex which they call a {\it Coxeter-like complex}, denoted $\Delta(G,S)$,
as follows.  The cells of $\Delta (G,S)$
are the cosets of parabolic subgroups of $G$, where a {\it parabolic 
subgroup} is defined to be any subgroup generated by a
subset of $S$.  A cell $gT$ is said to be in the closure of a cell $g'T'$ 
whenever $g'T' \subseteq gT$.  $\Delta (G,S)$ is the unique such  
regular cell complex with the 
property that each cell has the combinatorial type of
a simplex.  In other words, all lower
intervals in its face poset are Boolean algebras, implying that its face poset is
a simplicial poset (as in [St]), and $\Delta (G,S)$ is what is known as a cell complex of 
Boolean type (cf. [Bj]).  

Our focus will be on the case where the group $G$ is
the symmetric group $S_n$ and the minimal generating set 
$S$ is a set of $n-1$ transpositions that generate $S_n$. 
In this case, $\Delta(G,S)$ is a simplicial complex which is specified by the tree $T$ on
$n$ vertices 
having the edge set $\{ e_{i,j} |(i,j)\in S \}$, where $(i,j)$ denotes
the transposition swapping
$i$ and $j$.  Here we are using the fact that 
a set $S$ of transpositions generates $S_n$ if and only
if $\{ e_{i,j} | (i,j)\in S \} $ is a tree $T$ on $n$ vertices.  
Denote 
$\Delta(G,S)$ by $\Delta_T$  in this case, where $T$ is this associated tree. 
Observe  that the $i$-dimensional faces of $\Delta_T$ may be interpreted as 
the labelled forests obtained  by deleting $i+1 $ edges from $T$ and  
then assigning $|C|$ labels to each of the resulting components $C$ 
in such a way that each of the labels $1,\dots ,n$ is assigned to exactly one component.
When defining the ``inversions'' of a face later, 
we will think of such a labelled forest as a tree
whose vertices are the forest components and whose edges are the (deleted) edges
which connect the components, with the ``capacity'' of a component being the 
number of labels to be assigned to it. 

To allow for capacities larger than one, we study the more
general  type-selected Coxeter-like
complexes of [BR], denoted $\Delta_{(T,\bf{m})}$, where $T$ is a tree and $\bf{m}$ 
is a vector consisting of $|V(T)|$ nonnegative integers which specify the capacities
of the various vertices of $T$.  Denote the capacity of vertex $v$ by $\capac (v) $.
One way to define $\Delta_{(T,\bf{m} )}$ is as a Coxeter-like complex for the 
quotient group 
$S_{\sum_{v \in T} \capac (v)}  / S_{\capac (v_1)} \times \cdots \times 
S_{\capac (v_n)}$ with generating set $S$ consisting again of transpositions corresponding
to the edges of $T$.
Notice that one can also then make a 
combinatorial definition, analogous to the one given above for $\Delta_T $, again 
obtaining a face $F$ by deleting a set of edges $E^C(F)$ and then assigning 
the appropriate number of labels to  each of the resulting graph components, i.e. 
assigning exactly $\sum_{v\in C} \capac (v) $ labels to component $C$, where the 
labels are taken from the set $\{ 1,\dots ,\sum_{v\in T}  \capac (v) \} $.
Just as above, this is a simplicial complex with face containments of 
the form $\sigma \subseteq \tau $ if and only if $\sigma $ is obtained from $\tau $ by
merging neighboring components
Working in the generality of these
type-selected complexes $\Delta_{(T,\bf{m})}$ 
will be quite helpful to our analysis of the topological 
structure of skeleta of the complexes $\Delta_T$.  

Recall that the {\it chessboard complex} $M_{m,n}$ is the simplicial complex whose 
$i$-dimensional faces are the collections of $i+1$ mutually 
nonattacking rooks on an $m$ by $n$ chessboard.  
It was observed in [BR] that  $M_{m,n}$ is isomorphic 
to the Coxeter-like complex resulting from a 
tree with $m$ leaves,  each having capacity one and exactly one nonleaf 
vertex $v$ with $\capac (v) = n-m$.
The edges from $v$ to the leaves specify the $m$ 
rows of the chessboard while the $n$ labels specify the columns.  
A collection of non-attacking rooks  then corresponds to
an assignment of labels to the subset of the leaves whose incident edges have 
been chosen.

It was shown in [BLVZ] that $M_{m,n}$ is at
least $(\nu_{m,n} -1)$-connected.  Shareshian and Wachs later 
proved that $\tilde{H}_{\nu_{m,n} } (M_{m,n};\integers )\ne 0$ in [SW] and exhibited 3-torsion
in many cases.  Ziegler proved vertex decomposibility (and thus
shellability) of the  $\nu_{m,n}$-skeleta of chessboard complexes in [Zi] (see also [At]), 
while Friedman and Hanlon
determined $S_n$-module structure in 
[FH] and also showed there was no rational homology in degree $\nu_{m,n}$.  
See  [Wa] for a clear and quite comprehensive survey article regarding  
chessboard complexes and related complexes.

Our starting point was the following conjecture from [BR],
which we will prove in Theorem  ~\ref{proof-of-conjecture}.

\begin{conjecture}[Babson-Reiner]
\label{connectivity-bound}
If a tree $T$ has $n$ nodes and $b$ leaves, then
the Coxeter-like complex $\Delta_T$ is at least $(n-1-b)$-connected.
\end{conjecture}
In fact, we prove something stronger, namely 
that if $T$ has a collection $E $ 
of edges such that the tree $(T',\bf{m})$ associated to $(T,E)$
via Definition ~\ref{tree-assoc} 
satisfies $\capac (v) \ge deg(v)-1$ for each vertex
$v\in T$, and a similar condition holds for all subtrees,
then the $(|E|-1)$-skeleton of $\Delta_T$ is shellable, hence 
homotopy equivalent to a wedge of $(|E| -1)$-spheres.
The conjecture follows from the case in which $E$ is the set of
edges from parents to first children in a depth-first-search of a planar embedding
of $T$ with a leaf 
serving as the tree root.
In the course of proving this result, we also do several other things:

\begin{enumerate}
\item
axiomatize the notions of inversions and weak order for permutations
in a way that also makes sense for the 
labellings of any fixed tree (regarding permutations as  the labellings of a path)
\item
provide such an inversion function and weak order for any 
tree $(T,\bf{m})$ in which each 
vertex $v\in T$ satisfies $\capac(v) \ge deg(v) - 1$, i.e. any
so-called ``distributable capacity tree'' 
\item
prove that if $(T,\bf{m})$ admits an inversion function, then 
$\Delta_{(T,\bf{m})}$ is shellable
\item
extend a shelling constructed using (2) and (3) to one for the entire
$(n-b)$-skeleton of $\Delta_T$ for any tree $T$ with $n$ nodes and $b$ leaves 
\end{enumerate}
Our approach to shelling  via an inversion function also turns out to be 
related to the question of how to route packets of data efficiently 
on a network of computers which has a tree structure, viewing vertices as processors and edges as 
connections between them.  
Specifically, in Theorem  ~\ref{equiv-sort} we prove that  
the existence of an inversion function on the 
labellings of a tree 
is equivalent to that tree admitting greedy sorting 
based on this same inversion function; the requirements of an
inversion function ensure that
any series of local moves which swap labels at neighboring nodes that are 
out of order will terminate, and that all such series of local moves
terminate at the same fully-sorted tree.

Our setting is an idealized 
one which does not even begin to capture
the complexity of the networks of widest current interest such as the
internet graph, particularly since the graph of the 
internet is very far from being a tree.  However, one still might hope that our 
topological viewpoint on sorting and routing could offer some useful new insight. 
For instance, the connection we establish between greedy sorting and shellability
gives a topological obstruction to greedy sorting 
as follows: showing that $\Delta_{(T,\bf{m})}$
has nonvanishing homology in low degrees 
directly implies that greedy sorting 
based on an inversion set is not possible on the tree
$T$  with vertex capacities $\bf{m}$.  
Chessboard complexes already provide a class  of such trees, since
they are known to have nonvanishing homology in low enough
degrees to imply that their associated trees, namely stars, do not admit 
such greedy sorting algorithms
unless the central node of the star has capacity at least  its degree minus one.  
We refer readers to  [Le] for results in theoretical computer 
science regarding sorting and routing on various networks, 
though there are also many more recent results in this direction.

\section{Further terminology, notation and a key lemma}\label{notations}

The main focus of this paper is the Coxeter-like complex
$\Delta_T$ in which  $G$ is the symmetric
group $S_n$ and $S$ is a set of $n-1$ transpositions which generate 
$S_n$.  
The faces and incidences in this complex are as described in the introduction.
Given a face $F$, denote
by $E^C(F)$ the set of edges in $E(T)\setminus E(F)$, i.e. the edges deleted 
from $T$  to obtain the labelled forest representing $F$.  Thus, 
$\dim (F) = |E^C(F)|-1$.

\begin{definition}\label{tree-assoc}
Associate to a tree $(T,\bf{m})$ and a choice of subset $E$ of the set of  edges of $T$ 
a new tree $(T',\bf{m'})$ as follows.
The vertices of $T'$ are the connected 
components $C$ in the forest obtained by deleting from $T$ the edges in
$E$, and the edges of $T'$ are the edges in $E$.  The capacity of a 
component $C$ is the sum of its vertex capacities.
\end{definition}

\begin{figure}[h]
\begin{picture}(250,125)(40,10)
\psfig{figure=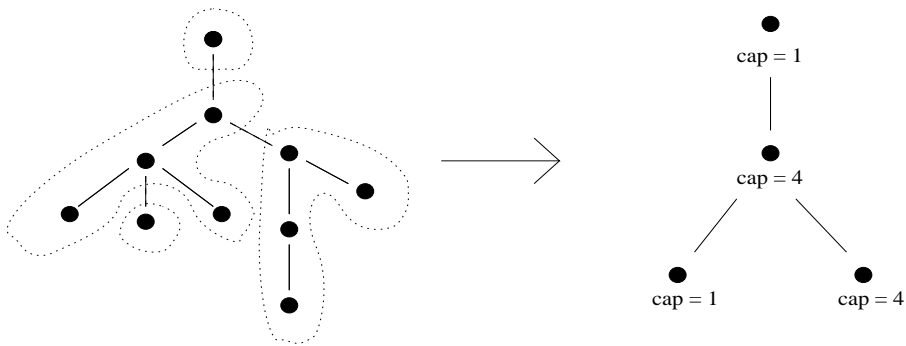,height=4.5cm,width=12cm}
\end{picture}
\caption{Tree derived from $(T,E)$}
\label{new-tree}
\end{figure}

Figure ~\ref{new-tree} gives one example of this process, with edges in $E$ depicted by dashed 
lines.  
For another  example, suppose that $T$ is a tree of four nodes, three of which are leaves, and
suppose each leaf has capacity one while the central node has capacity three.  Letting 
$E$ consist of any one edge of this tree, then the new tree $(T',\bf{m'})$ will be a path on two
nodes, one having capacity one and the other having capacity five.

\begin{definition}\label{distrib-def}
A tree $(T,\bf{m})$ has {\it distributable capacity} if each of its
vertices $v$ satisfies $\capac(v)\ge \deg(v) - 1$.
\end{definition}

For example, a tree with four nodes, three of which are leaves, has distributable 
capacity iff the non-leaf has capacity at least two.

\begin{definition}\label{face-assoc}
Given a tree $T$, a face $F$ in $\Delta_T$ is said to have {\it distributable capacity}
if the tree $(T',\bf{m'})$ obtained by applying Definition 
~\ref{tree-assoc} to tree $T$ and edge set $ E^C(F)$ has distributable capacity.
\end{definition}

\begin{lemma}\label{balanced-lemma}
The faces of 
$\Delta_T $ with distributable capacity
comprise a subcomplex.  Moreover, any component $C$ 
obtained by merging two neighboring components of a tree $T'$ with 
distributable capacity
satisfies $\capac (C) \ge \deg (C)$.
\end{lemma}

\proof
Proper faces  are obtained  by successively 
merging neighboring nodes.  If neighboring nodes 
$C_1,C_2$ are merged to form $C_1 \cup C_2$, and 
if $\deg(C_i) \le \capac (C_i) + 1$, for $i\in \{ 1,2\} $, then 
$C_1\cup C_2$ satisfies
\begin{eqnarray*}
\deg (C_1 \cup C_2) &=& \deg (C_1) + \deg(C_2) - 2 \\
&\le &  \capac (C_1) + \capac (C_2) \\
&<& \capac (C_1 \cup C_2) + 1.
\end{eqnarray*}
\EOP

Denote by $\overline{F}$ the simplicial complex comprised of  a face $F$ and all the 
proper faces of $F$. 
A {\it shelling} of a simplicial complex $\Delta $ is a
total order $F_1,\dots ,F_k $ on the facets of $\Delta $, i.e. on its maximal faces,  such that
for each $k>1$ the simplicial complex 
$\overline{F_k} \cap (\cup_{j < k} \overline{F_j} ) $ is a pure,
codimension one subcomplex of $\overline{F_k}$; recall that a simplicial complex 
is {\it pure} if all maximal faces are equidimensional.  
It is well-known (and
easy to show) that any shellable simplicial complex which is pure of dimension
$d$ for some $d\ge 2$ has homology concentrated in
top degree and is simply connected, hence is $(d-1)$-connected.
Recall that  the link of a face $G$ in a simplicial complex $\Delta $, denoted
$lk_{\Delta }(G)$, is the subcomplex
$lk_{\Delta }(G) = \{ F\in \Delta  | F\cap G = \emptyset, F\cup G \in \Delta  \} $, and 
that a shelling for $\Delta $ induces a shelling on the link of any face of $\Delta $.

Denote the $k$-skeleton of a simplicial complex $\Delta $ 
by $\Delta^{(k)}$.   We will prove that $\Delta_T^{(n-b)}$ has a pure, 
full-dimensional subcomplex  isomorphic to $\Delta_{(T',\bf{m})}$ for some $(T',\bf{m})$
with distributable capacity, and use the upcoming inversion function for $(T',\bf{m})$ to construct 
a shelling for $\Delta_{(T',\bf{m})}$ which extends to a shelling of $\Delta_T^{(n-b)}$.
By convention, we will always work in terms of a choice of 
rooted, planar embedding for $T$ which has a leaf as its root.
The {\it destination} of a label $i\in [n]$ in a tree $(T,\bf{m})$, 
denoted $dest(i)$, is the $a_i$-th 
vertex visited in a depth-first-search of $T$, where $\{ a_1,\dots ,
a_n \} $ is the integer-valued sequence with $a_1=1, a_i \le a_{i+1} 
\le a_i + 1$ for all $i$ which assigns 
$\capac (v) $ consecutive labels to
each $v$ as it is encountered in a depth-first-search.  Denote by $\lambda (T,\bf{m} )$
the set of labellings of a tree $(T,\bf{m})$ with labels $\{ 1,\dots ,\sum_{v\in T} \capac (v) \} $
which assign  $\capac (v) $ labels to vertex $v$ for each $v\in (T,\bf{m})$.
For further background in topological 
combinatorics and in Coxeter group theory, see [Bj3] and [Hu], 
respectively.

\section{Inversion function implies shelling}\label{implic-section}

In this section, we prove the following implications, after first defining what we mean by
inversion function:

\begin{enumerate}
\item
If $(T,\bf{m})$
admits an inversion function,  then $\Delta_{(T,\bf{m})}$ is shellable.
\item
If $\Delta_{(T,\bf{m})}$ is shellable via an inversion function for each 
$(T,\bf{m})$ having distributable capacity, 
then $\Delta_T^{(n-b)}$ is shellable for $T$ any tree with
$n$ nodes, $b$ of which are leaves.
\end{enumerate}

Theorem  ~\ref{no-full-shell} proves (1).  Theorem ~\ref{finish-shell} implies (2), but is 
more general.
In Section ~\ref{inversion-section}, we will construct an 
inversion function 
for all distributable capacity trees, enabling us to use (1) and (2) to give a  shelling for 
$\Delta_T^{(n-b)}$ and thereby deduce
the connectivity lower bound conjectured  in [BR].  

Generalizing Definition ~\ref{tree-assoc}, 
associate to any face $F$ in $\Delta_{(T,\bf{m})}$ a tree $T(F)$
by making a vertex for each 
connected component in $E(T)\setminus E^C(F) $ and an edge for each element of
$E^C(F)$, letting the capacity of any component be the sum of the capacities of the vertices
comprising that component.  Label the vertices of this tree with the collections of labels 
in the components of $E(T)\setminus E^C(F)$.
 Denote by $\lambda (T,\bf{m})$  the set of labellings of the tree 
$(T,\bf{m})$.  Now let $I$ be a function which assigns to each face
$F$ in $\Delta_{(T,\bf{m})}$ a collection of 
pairs of neighboring
components in $T(F)$, 
called the {\it inversion pairs} of $F$.  
Further require of $I$ that for each pair of neighboring components $C_1,C_2$,
there is a unique way to redistribute the labels collectively assigned to 
$C_1$ and $C_2$ so that $(C_1,C_2)$ is not an inversion 
pair of the resulting labelled tree $\tau $, i.e. so that $(C_1,C_2)\not\in I(\tau )$.  
This condition will in fact be subsumed by part 1 of 
Definition ~\ref{inv-fun-def}.

For neighboring components $i,j$ with $(i,j)\in I(\tau )$, 
define $(i,j)\tau $ to be the labelling obtained from $\tau $ by
redistributing the labels among $i$ and $j$ in the unique way so 
that $(i,j)\not\in I((i,j)\tau )$.  Make a covering relation
$(i,j)\tau \prec_{weak} \tau $ for each $(i,j)\in I(\tau )$ where $(i,j)$ are tree neighbors.

\begin{definition}\label{inv-fun-def}
A function $I$ as above is an {\it inversion function} for $(T,\bf{m})$
if $I$ has the following properties:
\begin{enumerate}
\item
Each face $F \in  \Delta_{(T,\bf{m})}$ is contained in a unique facet  
which is inversion-free on its restriction  
to each component of $F$.
\item
The transitive closure of 
$ \prec_{weak } $ is a partial order, which we then call the {\it weak order} with respect to $I$,
denoted $\le_{weak} $.
\item
These properties also hold on the restriction to any subforest.  
\end{enumerate}
\end{definition}

For example, suppose 
$T$ is a path with one end of the path  chosen as
tree root.  Say that a labelling of $T$ has an inversion pair $(v_i,v_j)$ exactly when 
$v_i$ is closer than $v_j$ to the root but some label assigned to $v_i$ is larger than some label
assigned to $v_j$.  It is easy to check that this meets the above requirements for an inversion
function, and in fact is equivalent to the usual notion of inversions between adjacent letters
in permutations.  However, it 
is much more difficult, and in many cases is impossible,  
to define an inversion function for other 
trees besides paths.   One of our main goals is to determine for which trees it is possible
to construct an inversion function, and one of the main results in this paper will be 
a sufficient  condition which will enable us to prove the conjecture of Babson and Reiner.

\begin{theorem}\label{no-full-shell}
If $(T,\bf{m}) $ admits an inversion function, then
$\Delta_{(T, \mm )} $ is shellable.
\end{theorem}

\begin{proof}
The facets in $\Delta_{(T,\mm )}$ may be viewed as 
labellings of $(T,\bf{m})$ 
with $\capac (v_i)$ labels assigned to the vertex $v_i$ for each $i$.  Denote by 
$F_w$ the facet associated to tree labelling $w$.

We first prove for each facet $F_w$ with $w$ not
the minimal element in  $\le_{weak }$
that $\overline{F}_w \cap (\cup_{u<_{weak} w} \overline{F}_u)$ is a pure
codimension one subcomplex of $\overline{F}_w$.
Implicitly here we use (2) to guarantee that $\le_{weak} $ is indeed a partial 
order.  By (1), merging neighboring components $j,j+1$ in $F_w$
for any  $(j,j+1)\in I(w)$ yields a 
codimension one face of $F_w$ also contained in the earlier facet
$ F_{(j,j+1) w} $.
Given any face $ G \in \overline{F}_w \cap
(\cup_{u<_{weak} w} \overline{F}_u )$, let $F_{u'}$ be a facet containing $G$
such that $F_{u'} <_{weak} F_w$.  Requirement (3) ensures that $u'$ may be obtained from
$w$ by a series of steps each eliminating an inversion between
two neighboring nodes, with the further restriction that these pairs
of neighboring nodes are each connected by edges
not belonging to $ E^C (G)$.  Thus, we obtain a codimension one face 
of $F_w$, denoted $G'$, with 
$G\subseteq G' \in \overline{F}_w \cap (\cup_{u<_{weak} w} \overline{F}_u )$ as follows:
merge two neighboring components of $F_w$
whose inversion may be eliminated as the first step in proceeding from 
$F_w$ to $F_{u'}$ in weak order.  

Now let $\le  $ be any linear extension of $\le_{weak} $.  
Consider $G$ any
face in $\overline{F}_n \cap (\cup_{m< n} \overline{F}_m )$.  
The earliest facet containing
$G$ is the unique one that is inversion-free on each component of $G$.
This will already come before $F_n$ in weak order, so
$\overline{F}_n \cap (\cup_{m<_{weak} n} \overline{F}_m ) = 
\overline{F}_n \cap (\cup_{m < n} \overline{F}_m )$ 
regardless of our choice of linear extension $\le $.
\end{proof}

\begin{remark}
It is easy to check that the case of a path amounts to exactly the shellability of the
type A Coxeter complex by using any linear extension of weak order to order the facets.
\end{remark}

It seems to be much easier for our upcoming
main example to confirm that a particular function $I$ is indeed an inversion function
by considering the much more 
detailed data of inversions between
non-neighboring nodes as well as neighbors.  In particular, the extra information
will make it 
easier to prove that the transitive closure of the set of covering relations is indeed a 
partial order.
Therefore, we make the following variation on the definition of inversion function, 
where now we let $I$ be a function assigning to each tree labelling a collection of 
(not necessarily adjacent) pairs of nodes.

\begin{definition}\label{inv-fun-def2}
A function $I$ assigning to each tree labelling a collection of pairs of nodes 
is a {\it  full inversion function} for $(T,\bf{m})$
if $I$ has the following properties:
\begin{enumerate}
\item
For each $\sigma\in\lambda (T,\bf{m}) $ and 
each $(i,k)\in I(\sigma )$, there is an inversion 
pair (denoted $(j,j+1)\in I(\sigma )$ to keep notation simple) 
between two neighboring components on 
the unique path from $i$ to $k$ in $(T,\bf{m})$.
\item
Each face $F$ is contained in a unique facet $F_m$
which is inversion-free on the restriction of $F_m$ to each 
component of $F$.
\item
The transitive closure of 
$ \prec_{weak } $ is a partial order, which we then call the {\it weak order} with respect to $I$,
denoted $\le_{weak} $.
\item
These properties also hold on the restriction to any subtree.  
\end{enumerate}
\end{definition}

Note that any full inversion function induces an inversion function by  assigning to 
each labelling its inversions which are between neighbors.

\subsection{Facet labelling giving rise to  $\Delta^{(n-b)}_T$ shelling}
\label{labelling-section}

\begin{figure}[h]
\begin{picture}(250,125)(-45,10)
\psfig{figure=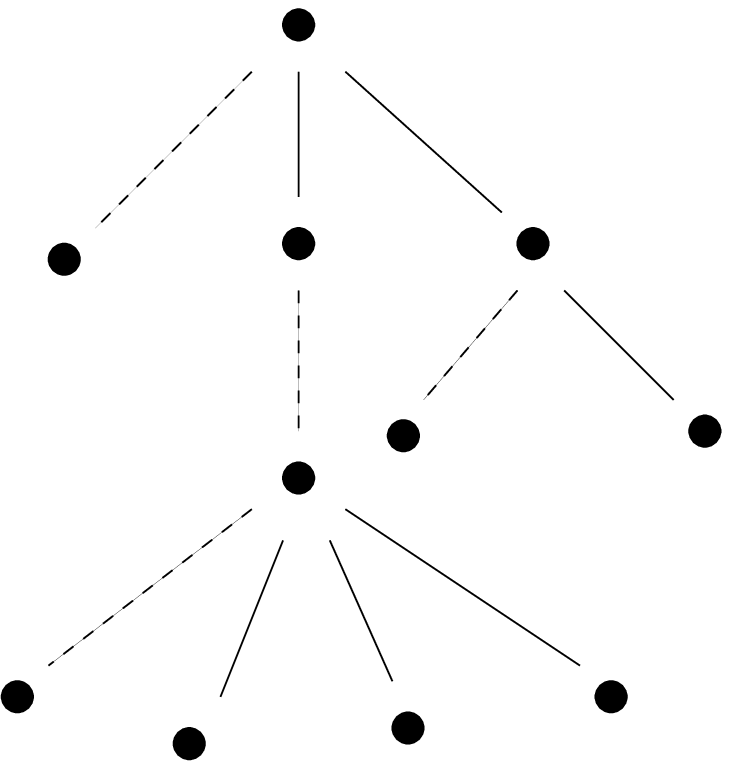,height=4.5cm,width=6cm}
\end{picture}
\caption{The set $E$ of edges to first children}
\label{first-children}
\end{figure}

Throughout this section,  $E  $ will always  
be a set of edges in a tree $T$ such that $(T,E)$ gives rise to a distributable 
capacity tree via Definition ~\ref{tree-assoc}.  Moreover, if we delete from $T$ any set $E'$ of
edges such that $E' \cap E = \emptyset $, and we
let $U$ be a connected component of the resulting
forest with $E(U)$ denoting the set of edges in $U$, 
then we also assume that
the tree $(T',\bf{m})$ obtained from $(U, E(U))$ via Definition ~\ref{tree-assoc} will 
also have distributable capacity.  Theorem ~\ref{finish-shell} will show  for such $(T,E)$
how to shell
$\Delta_T^{|E|-1}$.  Proposition ~\ref{c-0-case} will show that the set of edges from parents
to their first children in any depth first search of a tree whose root is a leaf will meet the above
conditions on a set $E$.  See Figure ~\ref{first-children} for what would be 
an example of such a set
$E$, except that in this example we have not chosen a leaf as the root .
We will use this special case 
to prove that $\Delta_T^{(n-b)}$ is shellable, for $T$ any tree with $n$ nodes and $b$
leaves.

Label each facet $F$ in $\Delta_T^{(|E|-1)} $
with coordinates $(c,S,I(F))$ as described next, and then order the facets by 
letting earlier coordinates take precedence over
later ones, using the orders specified below on individual coordinates.
Let $c= |E^C(F) \setminus E |$, and order this first coordinate linearly.
Let $S= E^C(F) \cap E$, and order this second coordinate by 
lexicographic order on the words obtained from the sets 
$S$ by listing the elements of $S$ in increasing order with respect to depth-first-search order
on the edges of the tree $T$.  

It will be convenient to first determine $c, S$ for a given facet $F$, 
and then describe the labels $I(F)$ for the various facets having this
fixed choice of $c$ and $S$.
In fact, we will need 
to fix somewhat more, namely the edge set $E^C(F)$, and then use the fact
that the facets of $\Delta_T^{(n-b)}$ with fixed edge 
set are the labellings of a fixed tree $(T',\bf{m})$ whose edges are exactly the edges
in $E^C(F)$.
Two facets $F_i,F_j$ may have 
the same coordinates $c$ and $S$ despite having $E^C(F_i)\ne E^C(F_j)$, but we deal with 
this by first ordering collections of facets using the labels $c$ and $S$, and then choosing any
ordering on the edge sets arising for a fixed $c$ and $S$ to extend to a total order on edge 
sets.  Next we describe how to assign sets $I(F)$ to the facets with a given edge
set  $E^C(F)$ and how to order these inversion sets.

When $c=0$, the tree $(T',\bf{m})$ obtained from $(T, E^C(F))$
by Definition ~\ref{tree-assoc} has
distributable capacity, as shown in Proposition ~\ref{c-0-case}.  Let $I(F)$ then be the set of
inversions given
by the inversion function for this distributable capacity tree (as provided in Section
~\ref{inversion-section}).  Partially order the facets  by the weak order
on their inversion sets, i.e., $\le_{weak}$,  and then choose
any total order extension of $\le_{weak}$ to obtain a total order on the
facets having a fixed edge set $E^C(F)$.    All total order extensions of $\le_{weak}$
are equivalent for purpose of shelling, since 
$\overline{F}_j \cap (\cup_{i<_{weak} j } \overline{F}_i ) = 
\overline{F}_j \cap (\cup_{i<j} \overline{F}_i)$ 
for each $j$ regardless of our choice of 
linear extension $<$, so there is no need to specify a choice.

For $c>0$, consider each component $U$ in the forest obtained from $T$ by deleting the
edges in $E^C(F)\setminus E$.  For each such $U$, consider the tree $(T',\bf{m})$ obtained from
$(U,E^C(F)|_U)$ via Definition ~\ref{tree-assoc}.  
This tree also has distributable capacity,
so let $I(F)$ be the union over all these components $U$
of the inversion sets for the various $U$.
Theorem ~\ref{finish-shell} will show that these sets $I(F)$ come
close enough to being inversion functions to allow us to extend 
the shelling of the subcomplex generated by the facets having $c=0$
to a shelling of the entire $(|E|-1)$-skeleton.

\subsection{Shelling skeleta of $\Delta_T$}

\begin{theorem}\label{finish-shell}
Suppose there exists $E \subseteq E(T)$
such that $(T,E)$ gives rise to a 
tree $(T',\bf{m})$ via Definition ~\ref{tree-assoc} which has distributable capacity.  Moreover,
for $U$ any subtree arising as a connected component in a forest obtained 
from $T$ by deleting any subset of the set of edges not in $E$, suppose that 
the tree $(T',\bf{m})$ obtained from $(U, E|_U)$ also has distributable capacity. 
Then $\Delta_T^{(|E|-1)}$ is shellable.
\end{theorem}

\begin{proof}
Let $F_1,\dots ,F_r$ be the facets of 
$\Delta_T^{(|E|-1)}$, ordered as described in Section ~\ref{labelling-section}.  We
will show this is a shelling in which the minimal face $G_m$ contributed by facet
$F_m$ consists of the following three types of edges in $E^C(F_m)$.
\begin{enumerate}
\item
 $E^C(F_m)\setminus E $
\item
$\{ \sigma \in E^C(F_m)\cap E | 
\sigma > \min (E \setminus E^C(F_m) ) \} $
\item
$\{ e_{i,j}\in E^C(F_m)\cap E  | (i,j) \in I(F_m) \} $
\end{enumerate}

The first thing to show is that deleting an edge of any of these three types from the 
facet $F_m$ yields a codimension one face that is also contained in an earlier 
facet.
Secondly, we must prove that
every face in $\overline{F}_m \cap (\cup_{m'<m} \overline{F}_{m'})$ is 
contained in
some such codimension one face.  Together, these will imply that
$\overline{F}_m \cap (\cup_{m'<m} \overline{F}_{m'})$ is a pure, codimension one 
subcomplex of $\overline{F_m}$, as is needed for the shelling.

To verify the first claim,
we show that omitting from $F_m$ any  vertex of $G_m$ will yield
a face $F$ that has codimension one in $F_m$ and
is shared with a facet $F_{m'}$ in which 
one of the coordinates in the facet labelling for $F_m$ has 
been decremented.  
If we delete from $F_m$ an edge in $E^C(F_m) \setminus E$, then there must be 
some edge in $E\setminus E^C(F_m)$ that may be added to $F$ to obtain
$F_{m'}$ with smaller $c$ coordinate and with $F\subsetneq F_{m'}$.  If  $F$ is
obtained by deleting from $F_m$ an edge from
$E^C(F_m)\cap E$ such that $E^C(F_m)$ lacks some earlier edge $e\in
E$, then we may insert $e$ into $F$ to obtain a facet with the same $c$ coordinate
as $F_m$ but a smaller $S$ coordinate.
Deleting an edge where there is an inversion between neighboring nodes means 
there will be an earlier facet in $\le_{weak} $ with this inversion eliminated, by 
Theorem ~\ref{inv-meets}.

Now we turn to the second claim.  
Suppose $G$ is a face in $\overline{F}_m\cap (\cup_{m'<m} 
\overline{F}_{m'})$
not contained in any  
codimension one face of $F_m$. 
Then $G$ must contain $G_m$ by (a).
In particular, $G$ must include $E^C(F_m)\setminus E$, hence
may only merge components of $F_m$ that are
connected by edges $e\in E$.
Moreover, $G$ may only merge components such that the edge $e$
connecting the pair of components
satisfies $e < \min (E\setminus E^C(F_m) )$.
Thus, $G$ is not contained in 
any $(|E|-1)$-faces in which either of the first two coordinates of $F_m$ have
been decremented from their value in $F_m$.  Additionally, $G$ is not
contained in any earlier facet $F_{m'}$ satisfying $E^C(F_m') \ne E^C(F_m)$ 
since $G$ contains all edges in $E^C(F_m)\setminus E$ as well as all edges 
in $E^C(F_m) \cap E$ that could possibly be replaced by earlier edges.
Finally, our first and third requirements for inversion sets
completely determine the distribution of labels in $F_m$ and force $F_m$
to come earlier than all other facets containing $G$, a contradiction to 
$G$ being shared with an earlier facet.  
\end{proof}

\begin{proposition}\label{c-0-case}
The tree $T$ together with the set $E$ of edges  from parents to first
children in a depth first search of $T$ with a leaf as root 
gives rise by Definition ~\ref{tree-assoc} to a tree $(T',\bf{m})$ with
distributable capacity.  Moreover, if we delete from $T$ any subset $E'$ of the edges  with
$E' \cap E = \emptyset $ and restrict 
to a component $U$ of the resulting forest, then 
the tree associated to $(U,U\cap E)$ by
Definition ~\ref{tree-assoc} also has
distributable capacity.
\end{proposition} 

\proof
Each component of the forest obtained by deleting these edges
has at least as many vertices as  it has 
edges from its vertices to their first children; the only other possible edge to another
component is from the root of the component to its parent, so $F$ has distributable 
capacity.   The same clearly holds on the restriction to any $U$, by virtue of our choice
of edge set $E$.
\EOP

\begin{corollary}\label{skeleton-shell-corollary}
If $T$ has $n$ nodes and $b$ leaves, then $\Delta_T^{(n-b)}$ is shellable.
\end{corollary}

\proof
There are $n-(b-1)$ edges from parents to first
children, since $T$ has a leaf as  root.
Letting $E$ be this collection of edges,
$(T,E)$ gives rise to a distributable-capacity
tree, and this property also restricts to subtrees as needed; this is
verified  in Proposition ~\ref{c-0-case}, so we may apply
Theorem ~\ref{finish-shell}.
\EOP

From this, the conjecture of [BR] is immediate:

\begin{theorem}\label{proof-of-conjecture}
If  $T$ has $n$ nodes and $b$ leaves, then
$\Delta_T$ is at least $(n-b-1)$-connected.
\end{theorem}

\section{Connection to greedy sorting on trees}

Define a {\it local sorting step} to be the unique redistribution of labels
between two neighboring nodes eliminating the  
inversion pair between the two nodes.  A tree labelling
is said to be {\it completely sorted} if it has no inversion pairs.

\begin{theorem}\label{equiv-sort}
Definition ~\ref{inv-fun-def}
may be rephrased in terms of greedy sorting algorithms as follows:
\begin{enumerate}
\item
Any tree
whose labels are not completely sorted
admits a local sorting step.
This is also true of the 
restriction to any subtree which is not completely sorted.
\item
All sequences of local sorting steps
lead to the same completely sorted labelling.  
This also holds for restrictions to subtrees, using only
local swaps within the subtree.
\item
Any series of local sorting steps will eventually terminate at a 
completely sorted tree.
\item
The above properties also hold for the restriction to any subtree.
\end{enumerate}
\end{theorem}

\proof
It is easy to see that the conditions of Definition ~\ref{inv-fun-def}
imply the conditions above, so we focus on proving the other direction.
The fact that any series of local sorting steps
terminates implies that the transitive closure of $\prec_{weak} $ is a partial 
order 
on the (finite) set of tree labellings, since otherwise there would
be a cycle enabling the
same series of sorting steps to be repeated indefinitely.
The fact that every series of
local swaps terminates with the same outcome implies that $\le_{weak }$
has a unique minimal element, giving the unique inversion-free labelling
of the entire tree.  The analogous requirement for restrictions
to subtrees implies Condition 2 also holds for all other faces.
\EOP

If one instead considers full inversion functions, i.e., Definition ~\ref{inv-fun-def2}, note then 
that condition 1in that definition is essentially equivalent to condition 1 of Theorem 
~\ref{equiv-sort}.
Combining Theorems ~\ref{equiv-sort} and  ~\ref{no-full-shell} yields:

\begin{corollary}\label{sorting-obstruction}
If $\Delta_{(T,\bf{m})}$ is not shellable, then $(T,\bf{m})$ does not 
admit a greedy sorting algorithm in the sense of  
Proposition ~\ref{equiv-sort}.
\end{corollary}

Thus, we obtain 
a homological obstruction to many  trees admitting greedy sorting 
algorithms based on inversion functions.  For instance, chessboard
complexes arise as the special case in which $T$ is a star, namely a
tree with a single vertex $v$ that is not a leaf.  
By Corollary ~\ref{sorting-obstruction},
known results on the homology of 
chessboard complexes from [SW] tell us
that stars with $\capac (v) < \deg (v) - 1 $ 
do not admit greedy sorting.

We should note that in the context of sorting/routing algorithms,
operations are typically performed simultaneously 
at the edges connecting many different pairs of processors.
However, an inversion function still indicates 
which such operations  would constitute 
progress in sorting data.  The distributable-capacity hypothesis
which will make possible an inversion function 
is very much in the spirit of results on routing in the sense that
allowing queues of bounded size to accumulate at nodes greatly improves
efficiency in algorithms.  

\section{An inversion function for distributable 
capacity trees}\label{inversion-section}

Let us begin with an example demonstrating some challenges to be overcome.
Trees are depicted with the root at the top  and
siblings ordered from left to right.
Edges in $E^C(F)$ are represented by dashed lines, while all other
edges in $E(T)$ are depicted by solid lines.   We say
$\lambda \in v_i$ when $\lambda $ is one of the labels assigned to vertex $v_i$.

\begin{example}\label{bottleneck-example}
A natural approach
would be to make an inversion $(i,j)$ 
whenever there are labels at positions $i$ and $j$ 
traversing the same edge in opposite directions on the most direct
routes to their destinations in depth-first-search order.
However, this inversion set is too sparse to satisfy
condition one of Definition ~\ref{inv-fun-def2} in general (and will also 
fail Definition ~\ref{inv-fun-def}).
\begin{figure}[h]
\begin{picture}(250,125)(-60,0)
\psfig{figure=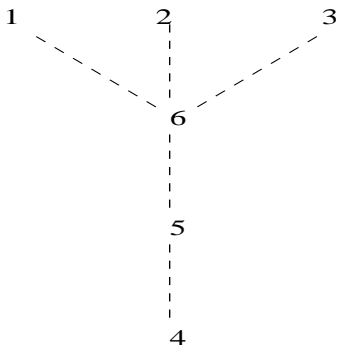,height=4.5cm,width=4.5cm}
\end{picture}
\caption{A queue of inverted labels}
\label{tree-queue}
\end{figure}
The example in Figure ~\ref{tree-queue} also illustrates a second issue, namely that of
contention.  There are inversion pairs between the vertices with label
pairs $(1,4),(2,5),(3,6)$, but it is not clear where there should 
be an inversion between neighboring nodes on 
the path between the locations of labels $1$ and $4$, as required by Definition
~\ref{inv-fun-def2}, condition one.  Notice also that this tree labelling would be 
inversion-free in the sense of Definition ~\ref{inv-fun-def}, but would not be the 
unique such labelling because the labelling which 
puts labels in exactly depth-first-search order would also be inversion-free.
\end{example}

The distributable capacity requirement will  handle this sort of contention
by ensuring that each 
vertex has enough capacity to allow distinct 
labels at the node to form inversions with its various neighbors, to exactly the 
extent that will be needed to define an inversion function;   Remark
~\ref{distrib-intuition} makes  this precise.  This simple 
idea will be crucial to the inversion function which is constructed 
over the remainder of this section. 
 
For any face $F\in \Delta_{(T,\bf{m})}$, 
define the  {\it components} of $F$ to be the components of the graph $G(F)$
obtained by deleting the edges in $E^C(F)$ from the tree 
$T$.  The inversions of $F$ will be the inversions of  
the labelled tree whose vertices are the 
components of $G(F)$ and whose edges are the edges in $E^C(F)$, letting the labels
of a vertex be the set of labels assigned to that component in $G(F)$.

\begin{remark}\label{distrib-intuition}
Let $T$ be a tree with distributable capacity and let $T'$ be any tree obtained from $T$ by
merging some neighboring nodes.  Then for any vertex $C\in T'$ comprised of more than
one vertex of $T$ and for any vertex $v\in T$
belonging to $C$, notice that the number of edges from $v$ to vertices of $T$ not in $C$
is at least as large as the capacity of $v$.   
\end{remark}

Call the edges $e_{v,w}$ from a fixed vertex $v$ 
the {\it capacity channels} of $v$.
Define the {\it path} of a label
$\lambda $ in a tree labelling to be the minimal path 
from $\lambda $'s position in the tree labelling 
to $dest(\lambda )$.  When we speak of an edge $e_{u,v}$ in the path 
of a label, by convention the directed path proceeds from $u$ to $v$.
A label $\lambda $ in a tree labelling
 {\it fills a capacity channel $e_{x,y}$} if the path of $\lambda $
includes the edge $e_{x,y}$.

The following notion 
will help us generalize ideas from linear arrays to trees, by 
specifying a hierarchy  of linear arrays within a tree:

\begin{definition}\label{coarsening-defn}
A {\it coarsening} with respect to a node $r$ (which we will regard as a ``local
root'')
is the decomposition of a rooted tree $T$ with respect to a marked node $r$ 
into as many as three parts as follows: (1) 
if  $r$ has multiple children, then 
there is a part $C_1$ consisting of
$r$'s first child  $v$ along with all descendants of $v$, (2) if $r$ has more
than two children, then there is a part $C_3$ 
consisting of $r$'s last child $w$ and all 
descendants of $w$, and 
(3) in any case there is a part $C_2$ consisting of 
$r$ and 
all remaining nodes of $T$ (some of which will not be descendants of $r$ if
$r$ is not the global root).
Now repeatedly subdivide $C_2$ in this fashion to 
obtain successive levels of coarsening on smaller and smaller subtrees, each of which 
includes $r$, until reaching a subtree in which $r$ has at most  one remaining child.
\end{definition}

\begin{definition}\label{path-dist-def}
The {\it coarsening-path-distance} of a label  $\lambda $ with 
respect to a particular coarsening at a chosen local root   $r$ 
is $|i-j|$, where $C_i$ and $C_j$ are the coarsening parts where the path of 
$\lambda $ begins and ends.
\end{definition}

While our upcoming inversion function may
seem  rather complicated, we are not aware of any simpler choice  that provably
works for more general  trees  than just paths.  
It would be interesting if a simpler inversion function exists in our level of generality,
i.e. for all distributable capacity trees.
 
 Part (1) of Definition ~\ref{distance-defn} and part 
 E1 of Definition ~\ref{label-inv-def} below are both forced by
 our upcoming base-labelling algorithm (which constructs the unique inversion-free labelling
 with proscribed labels).  Part (4) of Definition ~\ref{distance-defn} and parts V1-V4
 of Definition ~\ref{label-inv-def} are chosen so as to make inversion-free labellings tend towards
 ordering labels consistently with depth-first-seach order.  
  Other parts 
  constitute more arbitrary choices we have made  where some choice
 was needed.      
First we prioritize the capacity channels out of a vertex:

\begin{definition}\label{edge-priority}
Given a vertex  $u$, say that the edge $e_{u,w}$ is a {\it higher priority edge}
at $u$ than $e_{u,v}$ if either (1)
there is a coarsening with respect to local root $u$ such that
$u,v\in C_2$ but $w\not\in C_2$, or (2) there is a coarsening with respect to local
root $u$ such that  $w\in C_1$ and $v\in C_3$, or (3)
for all levels of coarsening we have $u,v,w\in C_2$ but with $v$ a child of $u$ and 
$w$ the parent of $v$.
\end{definition}

\begin{definition}\label{distance-defn}
Let $\mu ,\nu $ be two  labels in a labelled tree whose directed 
paths to their destinations intersect, 
letting $e_{x,y}$ be the final shared 
edge, when there is such an edge,  and otherwise letting  
$x$ be the unique shared vertex.
Then $\mu $ has
{\it higher priority} than $\nu $ at $x$,
or  in other words $\mu $ {\it travels farther} than $\nu $ from $x$,  
if any of the following conditions hold: 
\begin{enumerate}
\item
The path from $x$ to $dest(\mu )$ properly contains the path from $x$ to
$dest(\nu )$. 
\item
There is a shared edge $e_{x,y}$ as above, and
the coarsening-path-distance 
of  $\mu $ is greater than that of
$\nu $ with respect to some coarsening with $y$ as local root.
\item 
There is a shared edge $e_{x,y}$ as above, and the edge $e_{y,u}$ in the path
of  $\mu $ is a higher priority edge out of $y$ (cf. Definition ~\ref{edge-priority})
than the edge $e_{y,v}$ traversed by $\nu $.
\item 
There is a shared edge $e_{x,y}$ as above, 
$dest (\mu ) = dest(\nu )$, 
and either (a) $\mu < \nu $ with
$e_{x,y}$ proceeding from later to earlier position in
depth-first-search order, 
or (b) $\mu > \nu $
with $e_{x,y}$ proceeding from earlier to later position in 
depth-first-search order. 
\item
$\mu $ traverses a higher priority edge out of 
$x$ 
 than $\nu $ does.
\end{enumerate}
\end{definition}

Note that  in  Definition
~\ref{label-inv-def}, given next, 
two labels will never form a label inversion 
unless their paths share at least a vertex.  V1-V4 below will
deal with the various ways two paths 
may meet in just a vertex, while E1-E3 handle the ways two paths may 
share an edge or have the starting point of one path be contained in the other path.

\begin{definition}\label{label-inv-def}
A pair of  labels $(\lambda_i,\lambda_j)$ with $\lambda_i\in v_i$ and $ \lambda_j\in v_j$ 
forms a {\it label inversion pair} in a tree labelling 
$\sigma $, denoted $(\lambda_i,\lambda_j)\in I_{val}(\sigma )$, 
if any of the following conditions are met:
\begin{enumerate} [{E}1.]
\item 
$\lambda_i$ and $\lambda_j$ traverse the same edge in opposite 
directions on the paths from $v_i $ to $dest(\lambda_i ) $ and from
$v_j $ to $ dest(\lambda_j) $  \label{LA}
\item

$\lambda_j$ is on the path of $\lambda_i$, 
and   $\lambda_i$ fills a lexicographically smaller list
of unfilled capacity channels at $v_j$  (i.e. ones not filled
by higher priority labels according
to the label prioritization scheme described below) 
outward from $v_j$ than $\lambda_j$ does,
or there is a tie 
and $\lambda_i$ travels farther from $v_j$  than $\lambda_j$ does  \label{LB}. 
The lexicographic order is on  lists of indexing positions
for sublists of  the ordering on capacity channels given
within our label prioritization scheme below.
\item
The paths of $\lambda_i$ and $\lambda_j$ share at least one directed 
edge $e_{u,v}$  
for  $u\ne v_j$, and if  we let $e_{x,y}$ denote the final shared edge,   then 
$\lambda_j$ is closer than $\lambda_i$ to $y$ (in the sense of Definition
~\ref{distance-defn}), but
$\lambda_i$ travels farther than $\lambda_j$ does from $x$ (again as in Definition
~\ref{distance-defn}) \label{LC}
\end{enumerate}
Additionally, $(\lambda_i,\lambda_j)\in I_{val}(\sigma )$ if the
paths of $\lambda_i$ and $\lambda_j$ share a vertex $w$ but no
common edge, and
one of the following conditions is met:
\begin{enumerate}[{V}1.]
\item
$v_i,v_j$ are in parts $C_r,C_s$, respectively, of some 
coarsening with respect to a common ancestor $w$,
while $\lambda_i$ and $\lambda_j$ have destination parts 
$C_{d(r)}$ and $ C_{d(s)}$, respectively, with $r<s$ and $d(r)> d(s)$ \label{LE}
\item
$v_i \in C_r$ and $v_j\in C_s$ with $r<s$ for $C_r,C_s$ parts in some 
coarsening with respect to a common ancestor $w$,  
while $dest(\lambda_i),
dest(\lambda_j )$ are in the same part with respect to
all possible coarsenings, and $dest(\lambda_j) $ is an ancestor
of $w$ which is an ancestor of $dest (\lambda_i)$.
 \label{LF}
\item
$v_j$ is a descendent of $v_i$, $w$ is 
an ancestor of  $dest(\lambda_i)$ and 
$dest(\lambda_j)$, but $dest(\lambda_j )$ precedes $dest(\lambda_i )$ in 
depth-first-search order.  
\label{LG}
\item
$dest (\lambda_i) = dest (\lambda_j)$,
$\lambda_i < \lambda_j$, and $v_i$ comes later than $v_j$ in 
depth-first-search order \label{LH}
\end{enumerate}
\end{definition}

\begin{definition}\label{inv-def}
A node pair $(v_i,v_j)$ forms an {\it inversion pair} in a tree labeling $\sigma $, denoted
$(v_i,v_j)\in I(\sigma )$, if there are labels $\lambda_i\in v_i $ and
$  \lambda_j\in v_j$ forming a label inversion pair. 
\end{definition}

\begin{figure}[h]
\begin{picture}(250,125)(-15,10)
\psfig{figure=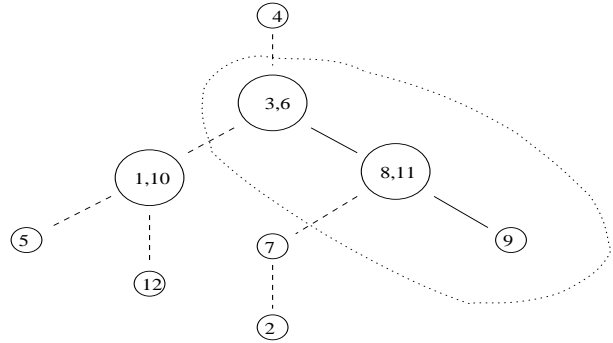,height=4.5cm,width=8cm}
\end{picture}
\caption{A labelling of a distributable capacity tree}
\label{tree-label}
\end{figure}

The tree labelling in 
Figure ~\ref{tree-label}  has inversions pairs $(v_1,v_2),
(v_2,v_3)$ and $ (v_6,v_7)$ between neighboring nodes resulting from label inversions 
of type E1, but has no inversions pairs in which both nodes belong to the set
$\{ v_2, v_6, v_9 \} $.

{\bf Label prioritization scheme:} 
Given a vertex $w$ 
and a collection of labels, first we prioritize edges, then  will
use the resulting ordering on
the edges  to prioritize labels.    Begin by 
ordering the edges  incident to $w$ according to Definition ~\ref{edge-priority}.
Then apply the list augmentation procedure below to add to the edge
list all remaining tree edges.
Now use this edge list to order the 
labels 
as follows.  Repeatedly give highest priority among labels not yet chosen 
to the label whose path from $w$ to its destination 
includes the lexicographically smallest
sublist of capacity channels  in the edge list that are 
not filled by any earlier  labels in the label list.  
Keep repeating to obtain an ordering on all the labels.

{\bf List augmentation procedure:} 
Proceed through a given ordered list of edges sequentially.  
Upon encountering an edge $e_{x,y}$,
append to the end of the list all  $e_{y,y'}$ for this fixed $y$ which are not
yet in the list, ordering those  travelling farthest from $y$  
(as in Definition ~\ref{distance-defn}, parts (2) and (3) applied to labels traversing
these edges)
earliest.  Continue  through the list  
until no further such augmentation is possible.

\begin{lemma}\label{unique-redistrib}
For each pair of vertices $v_i,v_j$ in a labelled tree, 
there is a unique redistribution of the set of labels collectively 
assigned to $v_i$ and $v_j$ such that $(v_i,v_j)$ is not an inversion pair
in the relabelled tree.
\end{lemma}

\begin{proof}
Let $P$ be the path  from $v_i$ to $v_j$.  For each capacity channel of $v_i$ or $v_j$
not involving any edges in $P$, the inversion-free
labeling will assign to $v_i$ (resp. $v_j$) the highest
priority label using that capacity channel, if any such label is available.  The remaining
labels must be assigned so as to avoid two labels traversing the same edge in 
opposite directions, which means that the vertex among $v_i,v_j$ having excess 
capacity will receive all the labels destined for it or outward from it across capacity
channels not involving any edges in $P$, as well as perhaps some additional labels.

For each label $\mu $ under consideration, there is a unique vertex $v\in P$ such that
the path from $v$ to $dest (\mu )$ only intersects $P$ in $v$.  This enables us to partition
the labels to be assigned to $v_i,v_j$ according to their associated vertices in $P$.
By Definition ~\ref{label-inv-def}, part E1, there is some vertex $x\in P$ such that all
labels associated to  $x' \in P$ with $x'$ closer than $x$ to $v_i$ must be assigned to 
$v_i$ in our inversion-free redistribution, and on the other hand all labels associated to 
$x''\in P$ having $x''$ closer than $x$ to $v_j$ must be assigned to $v_j$ in the 
inversion-free labelling.  
This just leaves the labels associated to $x$, but then
Conditions E3 and  V1-V4 of Definition ~\ref{label-inv-def}
determine  uniquely how to distribute to $v_i$ and $v_j$ 
these remaining labels. 
\end{proof}

\begin{figure}[h]
\begin{picture}(250,125)(-20,10)
\psfig{figure=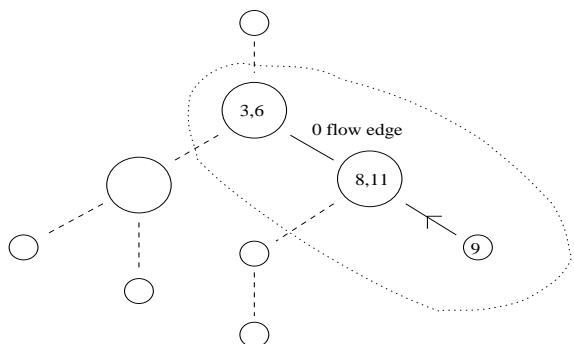,height=4.5cm,width=7.5cm}
\end{picture}
\caption{Base labelling of a component with a proscribed set of labels}
\label{base-label}
\end{figure}

The unique inversion-free labelling of a component with a proscribed set
of labels is its base-labelling, as defined next.  Figure ~\ref{base-label} gives an 
example of a base-labelling of a component.  Notice that this is a part of the tree
labelling given in Figure  ~\ref{tree-label}.

\begin{definition}\label{base-label-def}
The {\it base-labelling} of a component $C$ (of a forest obtained by deleting edges from a 
distributable capacity tree $T$) with a specified set of
labels is as follows.  If $C$ consists of a single vertex of $T$,
then all labels are assigned to 
the unique vertex of $C$, so henceforth assume $|C|>1$.  
Let $S$ be the maximal set of edges of $C$ which may be removed to yield a 
forest, each of whose components $C'$ has exactly as much capacity as the
number of labels destined either for 
$C'$ or across edges from $C'$ to components outside $C$.   Call this set of labels
$S(C')$.
Assign to each such $C'$ exactly the labels in $S(C')$, associating them to specific nodes 
within $C'$ as described next.

Obtain a directed graph on vertex set  $V(C')$
by orienting each edge of the induced subgraph on $V(C')$  from 
the portion of $C'$ with excess capacity
towards the portion with insufficient capacity.
Now we specify for each sink $v$ in this digraph
the set of labels to be
assigned to $v$; this will be a subset of the labels
either destined for $v$ or across edges $e_{v,w}$ for $w\not\in C'$.
For each capacity channel $e_{v,w}$ for $w\not\in C'$, 
assign to $v$ the highest priority label traversing $e_{v,w}$ 
according to  the  label prioritization scheme below. 
By Remark ~\ref{distrib-intuition}, we always have enough capacity to do this.
Then exhaust as much of the remaining capacity at $v$ as possible
by greedily assigning additional 
labels to $v$ whose paths involve edges from $v$ to other components, using  
the label prioritization scheme below to prioritize which labels are chosen. 
Break ties by choosing labels traveling farthest as in
Definition ~\ref{distance-defn}.   
Fill any remaining vacancies at $v$ with labels destined for $v$, doing so  in the 
unique way that is inversion-free according to Definition ~\ref{inv-def}.  Since $v$ is a 
sink, there will be enough such labels available to saturate $v$.

Once we thereby saturate $v$, remove $v$ 
from consideration, causing new nodes to become sinks and potentially splitting 
what remains of $C'$
into multiple components; keep repeating 
until each vertex in $C'$ has been assigned as many labels as its capacity.  
If $v\in C'$ has multiple edges in $C'$ leading to it just prior to its deletion, 
distribute as follows 
the unassigned labels  of $C'$ 
to the various components into which $C'$ splits upon  
deleting $v$.  Each component receives exactly as many of $v$'s 
excess labels as 
the capacity discrepancy across the edge from $v$ to that component;
specific excess labels  are 
distributed to these components in the unique way that is 
inversion-free according to  Definition ~\ref{inv-def}.   This is made possible by
the fact that 
each label is an excess label for at most one  sink at any given time -- because 
at any stage at
most one  sink is on the label's 
path to its destination within the entire tree (i.e. not just within $C'$).  
\end{definition}

The {\it base-labelling of a face}  is comprised of the base-labellings of its 
various components with their proscribed sets of labels.
The {\it inversions} of a face are the inversions in its base-labelling.

\begin{theorem}\label{inv-meets}
Definition ~\ref{inv-def} specifies  a full inversion function, hence its restriction to 
neighboring pairs of  nodes is an inversion
function.
\end{theorem}

\begin{proof}
It is straightforward 
to check the consistency and completeness
of Definitions ~\ref{base-label-def}, ~\ref{label-inv-def}, and 
that base-labellings are inversion free; this can be seen by drawing 
a few diagrams of different types of trees showing how 
Definition ~\ref{label-inv-def} covers the various ways the paths of two labels 
may intersect.   What remains is
to verify that Definition ~\ref{inv-def} satisfies the requirements from Definition
~\ref{inv-fun-def2} for a full inversion function.

The fact that an inversion $(v_i,v_k)$ forces an inversion
of two neighboring vertices on the path from $v_i$ to $v_k$
follows directly from the fact that 
traditional inversions in permutations have this property, 
together with our 
base-labelling algorithm  (i.e. Definition ~\ref{base-label-def})
which puts labels headed across the various 
capacity channels at nodes incident to these channels.  This gives Property (1)
of Definition ~\ref{inv-fun-def2}.
The first part of our base-labelling algorithm, i.e. the splitting into components $C'$ by 
deleting all possible edges having 0 flow, ensures that 
for each $G\subseteq F$, i.e. any face $G$ obtained by merging components of
$F$, that if $F$ is inversion free on the
restriction of $F$ to any component of $G$, 
then the base-labellings of $F$ and $G$ coincide on the restriction of $F$ to any
component of $G$.  This confirms 
Property (2) of 
Definition ~\ref{inv-fun-def2}.

Next, we  show that the transitive closure of $\prec_{weak}$ is 
a partial order, i.e. Property (3) of Definition ~\ref{inv-fun-def2}.  
This is accomplished  by verifying that the 
following integer-valued function $f$ on tree labellings satisfies
$f(\sigma ) < f(\tau )$ for each $\sigma \prec_{weak} \tau $.  For $\sigma $ a tree labelling, let
$N_{i,j}(\sigma )$ be the number of labels at $v_i$ in $\sigma $ that are assigned
to $v_j$ in the unique redistribution $\sigma '$
of labels between $v_i$ and $v_j$ satisfying $(v_i,v_j)\not\in I(\sigma ')$.  
Lemma ~\ref{unique-redistrib} guarantees the existence of such a $\sigma '$.  
Now let
$$f(\sigma ) = \sum_{(v_i,v_j)\in I(\sigma ) \atop i<j} 
N_{i,j}(\sigma )\cdot d(v_i,v_j),$$ 
with $d(v_i,v_j)$ 
denoting  distance in the usual graph-theoretic sense.  Let $d_{\sigma }(\lambda_i,\lambda_j) =
d(v_i,v_j)$ for $\lambda_i\in v_i$ and $\lambda_j\in v_j$ in the tree labelling $\sigma $.

Consider  $\sigma \prec_{weak} \tau $ for $\sigma = (r,r+1)\tau $.
If $\{ i ,j\} \cap \{ r,r+1 \} = 
\emptyset $, then $N_{i,j} (\tau ) = N_{i,j} (\sigma )$, so we need only
consider the contribution of pairs  $(i,j)$ with $\{ i,j\} \cap \{ r,r+1 \} \ne \emptyset $.
For  notational convenience, say
$i=r$ throughout the following argument; other cases are similar.

Suppose $j\ne r+1$.
Consider $(\lambda_{i_m},\lambda_{j_n}) \in I_{val}(\sigma )$ for 
$\lambda_{i_m} , 
\lambda_{j_n}$ located at $v_r,v_j$, respectively.  Then 
$(\lambda_{i_m},\lambda_{j_n}) \in I_{val}(\tau )$
with $d_{\tau }(\lambda_{i_m} ,\lambda_{j_n} ) \ge 
d_{\sigma }(\lambda_{i_m} ,\lambda_{j_n} )$ unless $(r,r+1)$ swaps 
$\lambda_{i_m} $ with some $\mu $, and we have 
$d_{\tau }(\lambda_{i_m},\lambda_{j_n}) + d_{\tau }(\mu,\lambda_{j_n})
= d_{\sigma }(\lambda_{i_m},\lambda_{j_n}) + d_{\sigma }(\mu ,\lambda_{j_n})$.
There may be a  choice of such $\mu $, but any 
choice will work, provided our collection of such choices comprises a matching 
on pairs of values in $I_{val}(\tau )$ being swapped  
between nodes $v_r$ and $v_{r+1}$. If $(\lambda_{i_m},\lambda_{j_n})\in I_{val}(\tau )$,
then $(\mu ,\lambda_{j_n} )\in I_{val}(\sigma )$ and $(\mu ,\lambda_{j_n} )\in
I_{val}(\tau )$.   If   
$(\lambda_{i_m},\lambda_{j_n})\not\in I_{val}(\tau )$, then
$\lambda_{i_m}$ is on the path from $\mu $ to $v_j$ in $\sigma $
with $(\mu ,\lambda_j)\not\in I_{val}(\sigma )$ and
$(\mu ,\lambda_j)\in I_{val}(\tau )$.    In each case, the pairs $(\lambda_{i_m},\lambda_{j_n})$
and $(\mu ,\lambda_{j_n})$ together contribute at least as much to $f(\tau )$ as to $f(\sigma )$.

Summing over all choices of $\{ i, j \} $ such that $\{ i,j\}  \ne \{ r, r+1 \} $
yields $f (\sigma ) - N_{r,r+1} (\sigma ) \le f(\tau ) - N_{r,r+1}
(\tau )$.  Finally,  
$N_{r,r+1}(\sigma )= 0$ while $N_{r,r+1}(\tau ) > 0$,
yielding $f(\sigma ) <  f(\tau )$, as desired.
Property (4) of Definition ~\ref{inv-fun-def2}, 
namely restriction to subtrees of the other properties, 
is immediate again from our definition of base-labelling.
\end{proof}

\section{Variations on the connectivity bound for $\Delta_T$}
\label{candidate-section}

A different choice of edge set $E$ which yields a 
capacity distributable tree will give a shelling of
$\Delta_T^{(j)}$ for 
$$
j =  \sum_{v\in T}\bigg\lfloor \frac{\deg (v)-1}{2} \bigg\rfloor .
$$
Namely, let $E$ consist of the edges from each vertex $v$ to its  
first $\lfloor \frac{\deg (v) -1}{2} \rfloor $ children, together with
the edge from the root to its only child.  Now apply Theorem
~\ref{finish-shell}.
In cases such as chessboard complexes
with nodes of high degree, this gives an
improved connectivity lower bound.  
Next we recover a shellability result from [Zi].

\begin{proposition}
The shelling for $\Delta_{(T,\bf{m})}$ when $(T,\bf{m})$ has
distributable capacity yields a shelling for $M_{m,n}$ for $n\ge
2m-1$.
\end{proposition}

\begin{proof}
Choose any edge set $E = \{ e_1,\dots ,e_k \}$ for $T$ a star with $m$
leaves.  Consider the subcomplex of $\Delta_{(T,\bf{m})} \cong M_{m,n}$
generated by those faces $F$ with 
$E^C (F) = \{ e_1,\dots ,e_k \} $.  These $F$ give rise 
to distributable
capacity trees if and only if $(m - k) + (n-m) \ge k - 1$.  Thus,
results of earlier sections yield shellability of 
$M_{m,n}^{(k-1)}$ for $k \le \frac{n+1}{2} $,  hence shellability 
of $M_{m,n}$ for $n\ge 2m-1$.
\end{proof}  

Our results do not  yield the optimal connectivity bound for 
chessboard complexes in general, regardless of our 
choice of $E$.   It is interesting to note that an $i$-dimensional 
type-selected complex $\Delta_{(T',\bf{m})}$ may not 
be shellable while the full $i$-skeleton of $\Delta_T$ may still be 
shellable.  It would be interesting to find the optimal connectivity 
bound for each $\Delta_T$, or more generally for each 
$\Delta_{(T, \bf{m})}$, the latter of which would include the chessboard 
complexes as a special case.

In  [He], we generalize
Ziegler's proof of vertex decomposability of skeleta to trees with one
non-leaf vertex as well as long exact sequences of Shareshian and Wachs from
[SW].  This enables a sharpness result, in the sense that we characterize
exactly which trees $(T,\bf{m})$ have shellable Coxeter-like complexes.
However, the question of a 
sharp connectivity bound remains open for general
Coxeter-like complexes.

Let $maxdeg(T)$ be the largest degree of any vertex in a tree $T$.  
Then the fact that Ziegler's vertex decomposition for
the $\nu_{m,n} $-skeleton of a chessboard complex generalizes to the
link of each face $F(v)$ together with the above theorem also suggests
the following question.

\begin{question}
Is the $(n-maxdeg(T))$-skeleton of $\Delta_T$ shellable? 
\end{question}   

Earlier sections prove a conjecture of Babson and Reiner by 
showing that the $(n-b)$-skeleton is shellable, where $b$ is the number
of leaves in $T$.  Notice that $b\ge maxdeg(T)$ holds, since 
each edge outward from a fixed vertex $v$ of maximal degree 
leads to a subtree with at least one leaf in it.  Thus, an affirmative
answer would give an improved connectivity bound for $\Delta_T$.

\end{document}